\documentclass[a4paper,12pt]{article}

\usepackage{amsmath}
\usepackage{times}
\usepackage{graphicx}
\usepackage[capitalize,nameinlink]{cleveref}
\usepackage{algorithm}
\usepackage{algorithmic}

\newcommand{\D}{\mathrm{d}}
\newcommand{\E}{\mathrm{e}}
\newcommand{\I}{\mathrm{i}}
\newcommand{\gnijk}[4]{g^{(#1)}_{#2,#3,#4}}
\newcommand{\gnbarijk}[4]{\overline{g}^{(#1)}_{#2,#3,#4}}

\crefname{section}{section}{sections}
\crefname{subsection}{subsection}{subsections}
\Crefname{section}{Section}{Sections}
\Crefname{subsection}{Subsection}{Subsections}

\Crefname{figure}{Figure}{Figures}

\crefformat{equation}{\textup{#2(#1)#3}}
\crefrangeformat{equation}{\textup{#3(#1)#4--#5(#2)#6}}
\crefmultiformat{equation}{\textup{#2(#1)#3}}{ and \textup{#2(#1)#3}}
{, \textup{#2(#1)#3}}{, and \textup{#2(#1)#3}}
\crefrangemultiformat{equation}{\textup{#3(#1)#4--#5(#2)#6}}%
{ and \textup{#3(#1)#4--#5(#2)#6}}{, \textup{#3(#1)#4--#5(#2)#6}}{, and \textup{#3(#1)#4--#5(#2)#6}}

\Crefformat{equation}{#2Equation~\textup{(#1)}#3}
\Crefrangeformat{equation}{Equations~\textup{#3(#1)#4--#5(#2)#6}}
\Crefmultiformat{equation}{Equations~\textup{#2(#1)#3}}{ and \textup{#2(#1)#3}}
{, \textup{#2(#1)#3}}{, and \textup{#2(#1)#3}}
\Crefrangemultiformat{equation}{Equations~\textup{#3(#1)#4--#5(#2)#6}}%
{ and \textup{#3(#1)#4--#5(#2)#6}}{, \textup{#3(#1)#4--#5(#2)#6}}{, and \textup{#3(#1)#4--#5(#2)#6}}

\crefdefaultlabelformat{#2\textup{#1}#3}

\title{A Fast Multipole Method for axisymmetric domains}

\author{Michael J. Carley}

\begin{document}

\maketitle

\begin{abstract}
  The Fast Multipole Method (FMM) for the Poisson equation is extended
  to the case of non-axisymmetric problems in an axisymmetric domain,
  described by cylindrical coordinates. The method is based on a
  Fourier decomposition of the source into a modal expansion and the
  evaluation of the corresponding modes of the field using a
  two-dimensional tree decomposition in the radial and axial
  coordinate. The field coefficients are evaluated using a modal
  Green's function which can be evaluated using well-known recursions
  for the Legendre function of the second kind, and whose derivatives
  can be found recursively using the Laplace equation in cylindrical
  coordinates. The principal difference between the cylindrical and
  Cartesian problems is the lack of translation invariance in the
  evaluation of local interactions, leading to an increase in
  computational effort for the axisymmetric domain. Results are
  presented for solution accuracy and convergence and for computation
  time compared to direct evaluation. The method is found to converge
  well, with ten digit accuracy being achieved for the test cases
  presented. Computation time is controlled by the balance between
  initialization and the evaluation of local interactions between
  source and field points, and is about two orders of magnitude less
  than that required for direct evaluation, depending on expansion
  order.
\end{abstract}

\section{Introduction}

Since its development~\cite{greengard-rokhlin87}, the Fast Multipole
Method (FMM) has become the algorithm of choice for a large class of
problems which can be expressed in terms of finding at a large number
of field points the potential generated by a large number of point
sources. This includes problems governed by the Poisson and Helmholtz
equations, including boundary integral problems in acoustics,
electromagnetism, and fluid dynamics, and volume integrals such as the
Biot--Savart integration which arises in electromagnetism and vortex
dynamics. The method is now highly developed with efficient
implementations available in two~\cite{ethridge-greengard01} and
three~\cite{langston-greengard-zorin11} dimensions for a range of
problems, using a variety of analytical tools for their formulation. 

Despite its power and importance, the FMM does not seem to have been
extended to non-Cartesian coordinate systems. In particular, there do
not seem to be formulations of the method which can be applied to
axisymmetric domains described using cylindrical coordinates. These
systems arise naturally in a range of applications and cylindrical
coordinates are a natural way to describe a system, or to apply
boundary conditions. The purpose of this paper is to present an
extension of the FMM for the Poisson equation to cylindrical
coordinates, motivated by applications in fluid dynamics, where
boundary and volume integral problems are governed by the Laplace
kernel. 

To the author's knowledge, there have been two previous studies which
are relevant to this problem. The first was the work of Strickland and
Amos~\cite{strickland-amos90,strickland-amos92} who developed an
accelerated method for the evaluation of the axisymmetric stream
function in vortex dynamics, equivalent to solving the Poisson problem
or evaluating the Biot-Savart integral. The authors used
single-precision arithmetic and a fifth order expansion of the Green's
function to achieve five-digit accuracy at a cost of~1--3\% of the
computational effort required for direct evaluation of the potential
and velocity fields. The authors do not seem to have extended their
method to the general case in a cylindrical domain.

More recently, in an unpublished thesis, Churchill~\cite{churchill16}
considered the general non-axisymmetric problem, motivated by the
analysis of boundary integrals on surfaces of revolution. He
identifies the particular difficulty in cylindrical coordinates, which
is that the Green's function is translation invariant in the axial
coordinate, but not in the radial, which complicates the evaluation of
interactions which arise in the FMM. He concludes that the problem can
be solved using a black-box~\cite{fong-darve09} or
generalized~\cite{gimbutas-rokhlin03} FMM, but does not present
results for the cylindrical problem. 

This paper presents a method for the fast evaluation of the potential
in a cylindrical domain, generated by a set of azimuthally varying
ring sources. The motivation is the Poisson or Biot--Savart problem in
fluid dynamics, with the expectation that the method is to be used in
boundary integral solvers~\cite[for example]{young-hao-martinsson12},
or in evaluating the velocity field due to a distribution of
vorticity~\cite{strickland-amos92}. The approach is essentially that
of a standard two-dimensional FMM, with modifications to the
evaluation of interactions to accommodate the lack of translation
invariance. It is assumed that any necessary Fourier transformation of
the inputs to the calculation has been performed, so that the starting
point is the location and radius of a set of circular sources, and the
coefficients of the Fourier series for a source distribution on each
circle. The main elements of the method are described briefly, with a
more detailed description of those parts which are particular to the
axisymmetric case, i.e.\ the evaluation of the Green's function and
its derivatives in cylindrical coordinates. Results are presented for
a problem with an increasing number of sources, to test the
performance of the method for convergence and speed.

\section{Analysis}
\label{sec:analysis}

\begin{figure}[htbp]
  \centering
  \includegraphics{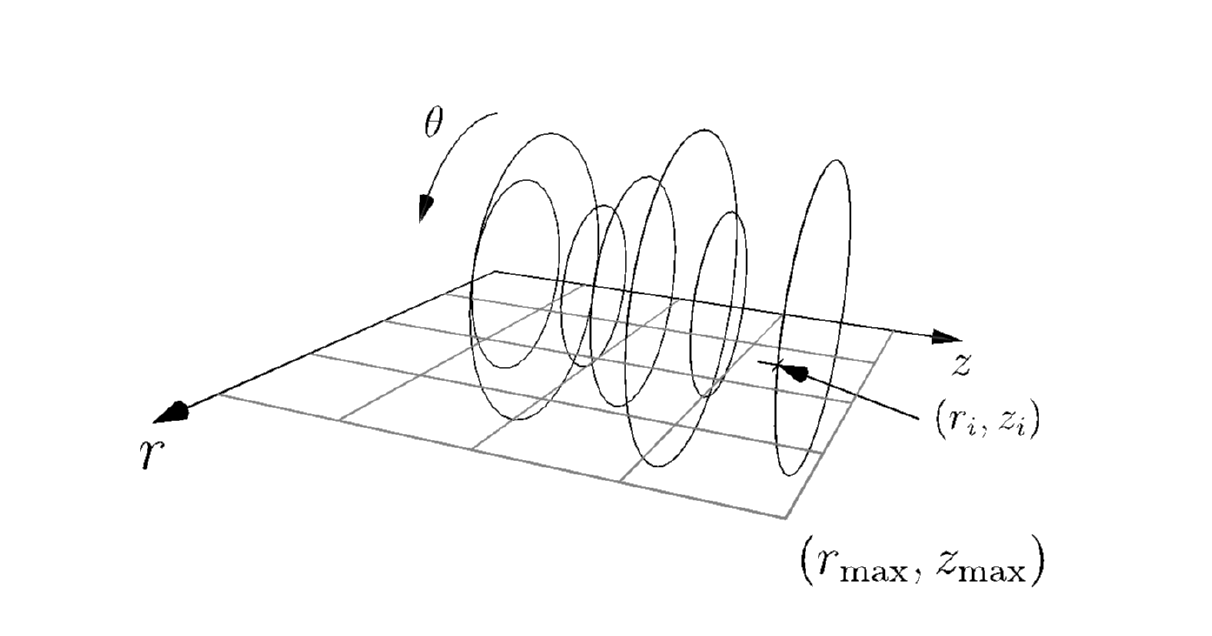}
  \caption{Basic problem: ring sources of radius $r_{i}$ are
    distributed at axial station $z_{i}$ in the domain
    $0\leq r\leq r_{\max}$, $0\leq z\leq z_{\max}$. Sources
    are given as the coefficients of the Fourier series of source
    strength as a function of $\theta$ at each $(r_{i},z_{i})$}
  \label{fig:analysis:problem}
\end{figure}

The problem to be solved is sketched in
\cref{fig:analysis:problem}. We assume that any necessary
preprocessing has been performed to reduce the system to the form
shown here. A set of ring sources with a common axis are distributed
throughout the domain $(r,z)$. The source strength $s(\theta_{1})$
on a ring at location $(r_{1},z_{1})$ is given by the Fourier series,
\begin{align}
  \label{equ:analysis:source}
  s(\theta_{1})
  &=
  \sum_{n=-N}^{N}
  S^{(n)}\E^{\I n \theta_{1}},
\end{align}
where subscript~$1$ denotes source coordinates. The potential
$\phi(r,z)$ due to the source is given by integration over
$\theta_{1}$,
\begin{align}
  \phi(r,\theta,z)
  &=
  \int_{0}^{2\pi} \frac{s(\theta_{1})}{4\pi R}\,\D\theta_{1},\\
  R^{2} &= r^{2} + r_{i}^{2} - 2rr_{1}\cos(\theta-\theta_{1}) +
  (z-z_{1})^{2}.\nonumber
\end{align}
Under an elementary transformation, and substituting the Fourier
series for $s(\theta_{1})$,
\begin{align}
  \phi(r,\theta,z)
  &=
  \sum_{n=-N}^{N}
  S^{(n)}\E^{-\I n\theta}
  \int_{0}^{2\pi} \frac{\E^{\I n\theta_{1}}}{4\pi R}\,\D\theta_{1},\\
  R^{2} &= r^{2} + r_{1}^{2} - 2rr_{1}\cos\theta_{1} +
  (z-z_{1})^{2}.\nonumber  
\end{align}
The potential field $\phi(r,z)$ can then be expressed as a Fourier
series in $\theta$,
\begin{align}
  \label{equ:analysis:phi}
  \phi(r,\theta,z)
  &=
  \sum_{n=-N}^{N}
  \Phi^{(n)}(r,z)\E^{\I n \theta},\\
  \Phi^{(n)} &= S^{(n)}G^{(n)}(r,r_{1},z-z_{1}),\\
  G^{(n)}(r,r_{1},x)
  &=
  \int_{0}^{2\pi} \frac{\E^{\I n\theta_{1}}}{4\pi R}\,\D\theta_{1}
  =
  \int_{0}^{2\pi} \frac{\cos n\theta_{1}}{4\pi R}\,\D\theta_{1},\\
  R^{2} &= r^{2} + r_{1}^{2} - 2rr_{1}\cos\theta_{1} + x^{2}.\nonumber
\end{align}
The source term $s(\theta_{1})$ is assumed real, so that the complex
Fourier coefficients are related by
$S^{(-n)}=\left(S^{(n)}\right)^{*}$. Henceforth, all computations will
be performed for $n\geq 0$ and the conjugate relationship will be
assumed.

We refer to $G^{(n)}(r,r_{1},x)$ as the \emph{modal Green's function}
relating Fourier coefficients of the source to those of the potential
at some other point. \cref{sec:gfunc} gives details of the evaluation
of $G^{(n)}(r,r_{1},x)$ and of its derivatives. In particular, it is
shown that it can be expressed exactly as
\begin{align}
  \label{equ:analysis:gfunc}
  G^{(n)}(r,r_{1},x)
  &=
  \frac{Q_{n-1/2}(\chi)}{2\pi\sqrt{rr_{1}}},  
\end{align}
where $Q_{\nu}(\chi)$ is a Legendre function of the second
kind~\cite{cohl-tohline99}.

Given a distribution of sources at locations $(r_{i},z_{i})$,
$i=1,\ldots,N_{s}$, each of which has a set of Fourier coefficients
$S_{i}^{(n)}$, $n=0,\ldots,N$, the problem to be solved is then the
approximate evaluation of the sum
\begin{align}
  \label{equ:analysis:sum}
  \Phi^{(n)}(r_{j},z_{j})
  &=
  \sum_{i=1}^{N_{s}}
  S^{(n)}_{i}G^{(n)}(r_{j},r_{i},z_{j}-z_{i}),
\end{align}
for field points $(r_{j},z_{j})$, $j=1,\ldots,N_{f}$. 

\subsection{Outline of the FMM}
\label{sec:analysis:outline}

\begin{figure}[htbp]
  \centering
  \includegraphics{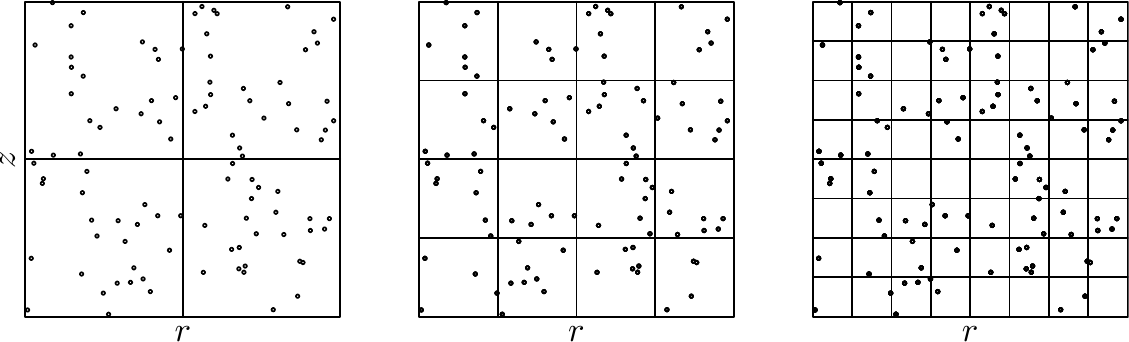}
  \caption{Recursive subdivision of the source domain into $4^{\ell}$
    boxes for $\ell=1,2,3$, left to right}
  \label{fig:fmm:division}
\end{figure}

The Fast Multipole Method is well established and there are numerous
guides to its algorithm and implementation. Here we give an outline of
the method in order to present the necessary terminology and to
indicate those parts of the algorithm of this paper which differ from
existing methods. From the previous section, we recall that the
objective is to approximately evaluate the sum
\begin{align*}
  \Phi^{(n)}(r_{j},z_{j}) &= \sum_{i=1}^{N_{s}}
  S_{i}^{(n)}G^{(n)}(r_{j},r_{i},z_{j}-z_{i}),\quad n=0,\ldots,N,
  \quad j=1,\ldots,N_{f},
\end{align*}
given a list of source and field points $(r_{i},z_{i}$ and
$(r_{j},z_{j})$ respectively. The first operation of the FMM is the
sorting of points and their representation in a quadtree data
structure, formed by repeated subdivision of the domain
$0\leq r_{i,j}\leq r_{\max}$, $0\leq z_{i,j}\leq z_{\max}$,
$z_{\max}=r_{\max}$. \cref{fig:fmm:division} shows the repeated
halving of intervals. At each \emph{level} of subdivision $\ell$, the
domain is divided into $2^{\ell}\times2^{\ell}$ boxes, to a maximum
level called the \emph{depth} of the tree $d$. Boxes in the tree can
be indexed by their location $(i,j)$ on the grid at a given level, or
by their Morton index. 

\begin{figure}[htbp]
  \centering
  \includegraphics{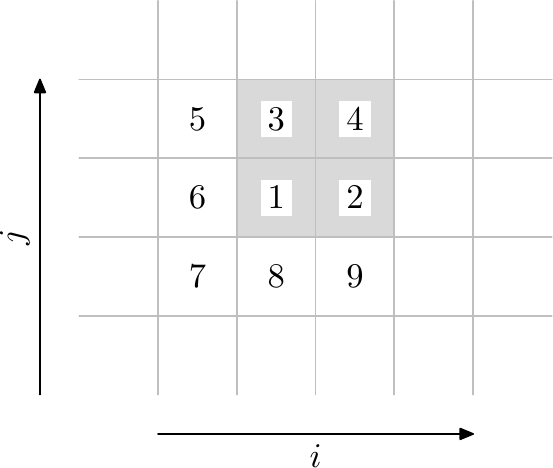}
  \caption{Terminology for relationships between boxes: boxes~1--4 are
  children of the shaded box; the shaded box is the parent of
  boxes~1--4; boxes~1--9 are neighbors of box~1}
  \label{fig:fmm:boxes}
\end{figure}

\cref{fig:fmm:boxes} gives the terminology for relationships between
boxes. A box at grid location $(i,j)$ at level $\ell$ has four
\emph{child} boxes at grid locations $(2i,2j)$, $(2i,2j+1)$,
$(2i+1,2j)$, $(2i+1,2j+1)$ at level $\ell+1$, and has a \emph{parent}
box at grid location $(i/2,j/2)$ at level $\ell-1$. A box at the
finest level of subdivision, depth $d$, has no children and is called
a \emph{leaf} box. 

\emph{Neighbors} of box $(i,j)$ are boxes at the same level $\ell$
which share at least a vertex with the box, including the box
itself. In \cref{fig:fmm:boxes}, boxes~1--9 are neighbors of
box~1. All other boxes are said to be in the \emph{far field} of
box~1. The basic principle of the FMM is to separate far-field from
neighbor interactions and evaluate the far-field terms in any box
using an accelerated summation. 

This is achieved in the first instance by evaluating the field due to
sources inside a box using an approximate expansion which is faster
than direct evaluation of the sum
\begin{align*}
  \Phi^{(n)}(r,z)
  &=
  \sum_{i}S_{i}^{(n)}G^{(n)}(r,r_{i},z-z_{i}),
\end{align*}
where the summation is taken over all sources inside a box. At field
points sufficiently far from the box, the Green's function for a
source is well approximated by its Taylor series, truncated to some
order $M$,
\begin{align}
  \label{equ:analysis:expansion}
  G^{(n)}(r+\Delta r,r_{1}+\Delta r_{1},x+\Delta x)
  &\approx
  \sum_{m=0}^{M}
  \sum_{i+j+k=m}
  \gnbarijk{n}{i}{j}{k}(\Delta r)^{i}(\Delta r_{1})^{j}(\Delta
  x)^{k},\\
  \gnbarijk{n}{i}{j}{k}
  &=
  \frac{1}{i!}\frac{1}{j!}\frac{1}{k!}
  \frac{\partial^{i+j+k}}{\partial r^{i}\partial r_{1}^{j}\partial
    x^{k}} G^{(n)}(r,r_{1},x).\nonumber    
\end{align}
Details of the evaluation of the derivatives of $G^{(n)}$ are given in
\cref{sec:gfunc}.

The field at a point $(r,z)$ due to a source at
$(r_{1}+\Delta r_{1}, z_{1}+\Delta z_{1})$ in a box with center
$(r_{1},z_{1})$ is then given by
\begin{align}
  \label{equ:field:expansion}
  \Phi^{(n)}(r,z)
  &\approx
  S^{(n)}
  \sum_{m=0}^{M}
  \sum_{i+j=m}
  (-1)^{j}
  \gnbarijk{n}{0}{i}{j}(\Delta r_{1})^{i}(\Delta
  z_{1})^{j},
\end{align}
noting that $x=z-z_{1}$ and $\partial/\partial
z_{1}=-\partial/\partial x$. Summing over all sources contained in the
box,
\begin{align}
  \Phi^{(n)}(r,z)
  &\approx
  \sum_{m=0}^{M}
  \sum_{i+j=m}
  (-1)^{j}
  \gnbarijk{n}{0}{i}{j}
  S^{(n)}_{ij},
\end{align}
where the moments $S_{ij}$ are given by
\begin{align}
  \label{equ:analysis:moments}
  S^{(n)}_{ij}
  &=
  \sum_{q} S_{q}^{(n)} (\Delta r_{q})^{i}(\Delta z_{q})^{j}.
\end{align}

\begin{figure}[htbp]
  \centering
  \includegraphics{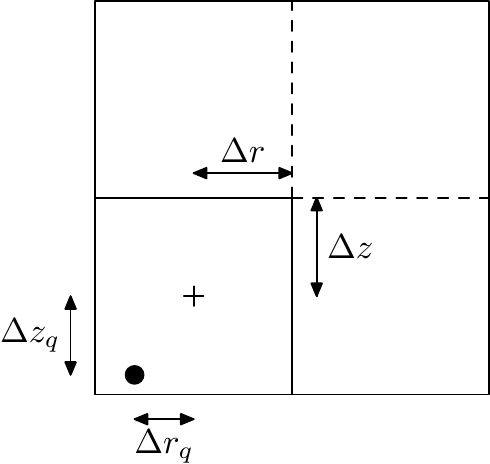}
  \caption{Evaluation and shifting of source moments: moments about
    the center of the leaf box are computed using
    \cref{equ:analysis:moments}; those about the center of the parent
    box are computed by shifting the child box moments by $(\Delta
    r,\Delta z)$ using \cref{equ:analysis:shift}.}
  \label{fig:analysis:moments}
\end{figure}

To initialize the source data in the first stage of the FMM, the
moments $S_{ij}^{(n)}$ are computed for each leaf box in the tree. In
the \emph{upward pass}, the moments at boxes in each level $\ell$ are
computed, for $\ell=d-1,\ldots,1$. This can be achieved without
requiring direct evaluation of moments from source data, by combining
moments from child boxes to generate moments in their parent box,
\cref{fig:analysis:moments}. Moments about the center of a child box
at displacement $(\Delta r, \Delta z)$ contribute to the moments about
the center of their parent box via
\begin{align}
  \label{equ:analysis:shift}
  S_{ij}
  &=
  \sum_{q=0}^{i}
  \sum_{u=0}^{j}
  \binom{i}{q}
  \binom{j}{u}
  (-\Delta r)^{q}
  (-\Delta z)^{u}
  S_{i-q,j-u}',
\end{align}
where the superscript $(n)$ has been suppressed for clarity. When the
upward pass has been completed, each box at levels $\ell=1,\ldots,d$
contain a set of moments which can be used to estimate the potential
in the far field of the box.

In the next stage of the FMM, the \emph{downward pass}, each box is
assigned a \emph{local expansion} which can be used to evaluate the
potential inside the box due to sources which lie in its far
field. The core of the FMM is the use of the most efficient expansion
possible at any level to evaluate the far-field terms in any box. The
field in a box centered at $(r,z)$ is given by
\begin{align}
  \label{equ:analysis:local}
  \Phi^{(n)}(r+\Delta r,z+\Delta z)
  &=
  \sum_{m=0}^{M}
  \sum_{k+\ell=m}
  \Phi_{k\ell}^{(n)}
  (\Delta r)^{k}(\Delta z)^{\ell},
\end{align}
where the expansion coefficients $\Phi_{k\ell}^{(n)}$ are evaluated
from the contributions of sources in boxes which interact with the
field box. The order $M$ of the local expansion at any level is not
required to be the same as the order of the source expansions, but has
been set equal for the calculations presented in this paper.

\begin{figure}[htbp]
  \centering
  \includegraphics{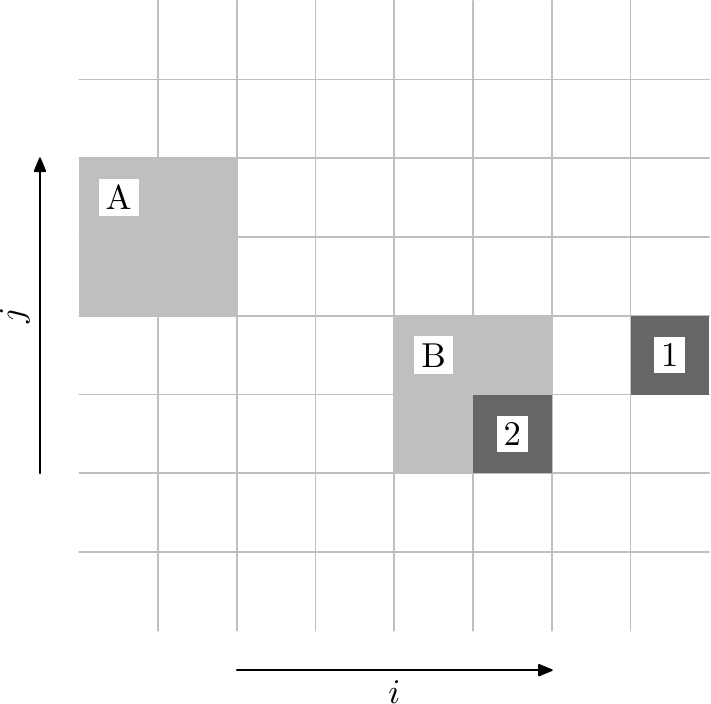}
  \caption{Generation of a local field in a box: the local expansion
    of the field in box B due to sources in box A is computed using
    the S2L operation, \cref{equ:analysis:s2l}; the local expansion in
    box~2 is found by shifting the expansion in box B to box~2 and
    adding the contribution from box~1, evaluated using an S2L
    operation.}
  \label{fig:analysis:local}
\end{figure}

\cref{fig:analysis:local} shows the main operations involved, for the
evaluation of a local expansion in box~2, which has parent box B. On
the downward pass, the local expansion in box B is found by adding the
contribution from sources in boxes which are well separated from B,
such as box A. This contribution is found using the
\emph{shift-to-local} or S2L operation. The local expansion in B is
used to generate the local expansion in each of its child boxes,
including box~2. Box~2 then has its local expansion incremented by the
contribution of boxes with which it interacts, such as box~1. At the
end of the downward pass, each leaf box has a local expansion which
accounts for the contribution of all sources lying outside its
neighbors. 

The two operations to be implemented here are the S2L and the
parent-to-child shift of the local expansion. The S2L shift is
found by differentiating \cref{equ:field:expansion},
\begin{align}
  \label{equ:analysis:s2l}
  \Phi^{(n)}_{k\ell}
  &=
  \sum_{i,j}
  (-1)^{j}
  \binom{j+\ell}{j}
  S^{(n)}_{ij}\gnbarijk{n}{k}{i}{j+\ell},
\end{align}
which can be implemented as a BLAS level~2 operation

The local expansion in a child box at displacement $(\Delta r,\Delta
z)$ is given from the parent box expansion by,
\begin{align}
  \label{equ:analysis:l2l}
  \Phi_{ij}
  &=
  \sum_{q=0}
  \sum_{u=0}
  \binom{i+q}{q}
  \binom{j+u}{u}
  (-\Delta r)^{q}
  (-\Delta z)^{u}
  \Phi_{i+q,j+u}',
\end{align}
where terms $\Phi_{i,j}'$ are coefficients of the parent box local
expansion. 

\subsection{Evaluation of interactions}
\label{sec:analysis:interactions}

The outline of the Fast Multipole Method presented in
\cref{sec:analysis:outline} contains the main elements of a generic
FMM which are familiar from existing implementations. In this section,
we describe the part of the algorithm which is particular to the
cylindrical domain, the S2L operation for the modal Green's
function. In existing, Cartesian, methods, the translation operators
are invariant with respect to shifts in the coordinate system. As
noted by Churchill~\cite{churchill16}, however, this is not true for
the modal Green's function $G^{(n)}$, which is invariant for shifts in
the axial coordinate $z$ but not for displacements in radius $r$. This
increases the number of orientations for which shift operators must be
computed, though there are still some symmetries which can be
exploited to reduce the workload. 

\begin{figure}[htbp]
  \centering
  \includegraphics{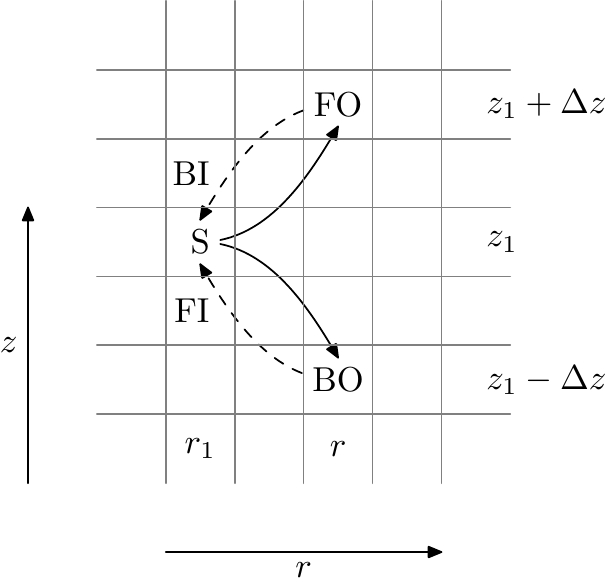}
  \caption{Reuse of Green's function derivatives for forward/backward
    and inward/outward interactions. The source box S generates local
    expansions at $r>r_{1}$ through a forward-outward (FO) and a
    backward-outward (BO) operation, indicated by solid arrows. The
    boxes at greater radius generate local expansions in box S via
    backward-inward (BI) and forward-inward (FI) operations. All four
    operations are computed from the same expansion of the Green's
    function.}
  \label{fig:analysis:s2l}
\end{figure}

Recall that the modal Green's function is given by
\begin{align*}
  G^{(n)}(r,r_{1},x)
  &=
  \frac{Q_{n-1/2}(\chi)}{2\pi\sqrt{rr_{1}}}, \,
  \chi = \frac{r^{2}+r_{1}^{2}+x^{2}}{2rr_{1}},\,x = z-z_{1}.
\end{align*}
This is invariant under translations in $z$ and is symmetric in $r$
and $r_{1}$, a fact which is exploited in the recursion relations for
derivatives in source and field
coordinates~\cite{strickland-amos90,strickland-amos92}. To take
advantage of this symmetry, we introduce some terminology to describe
translation operations. If we assume that for S2L operations, we
evaluate derivatives of $G^{(n)}(r,r_{1},x)$ for the $r\geq r_{1}$,
$x\geq 0$. Then, from \cref{fig:analysis:s2l}, we can derive
translation operators for four different cases. These correspond to
shifts in the positive or negative (forward or backward) axial
direction, and from greater to smaller radius (outward or inward).

The basic operator, which uses the derivatives of the Green's function
without modification, is the \emph{forward-outward} or FO shift,
\cref{equ:analysis:s2l}. 
\begin{align*}  
  \Phi^{(n)}_{k\ell}
  &=
  \sum_{i,j}
  (-1)^{j}
  \binom{j+\ell}{j}
  S^{(n)}_{ij}\gnbarijk{n}{k}{i}{j+\ell}.
\end{align*}
For the \emph{backward-outward} or BO shift, $x<0$,
$\partial/\partial z_{1}=\partial/\partial x$
and
$\partial/\partial z=-\partial/\partial x$, yielding
\begin{align*}  
  \Phi^{(n)}_{k\ell}
  &=
  \sum_{i,j}
  (-1)^{\ell}
  \binom{j+\ell}{j}
  S^{(n)}_{ij}\gnbarijk{n}{k}{i}{j+\ell}.
\end{align*}

To evaluate the inward shifts, we exchange $r$ and $r_{1}$ and swap
the corresponding indices. For the \emph{forward-inward} S2L
operation,
\begin{align*}  
  \Phi^{(n)}_{k\ell}
  &=
  \sum_{i,j}
  (-1)^{j}
  \binom{j+\ell}{j}
  S^{(n)}_{ij}\gnbarijk{n}{i}{k}{j+\ell},
\end{align*}
which gives the contribution of the sources in the box at larger
radius $r$ to the local expansion in the box at smaller radius
$r_{1}$. Finally, the \emph{backward-inward} operator is given by
\begin{align*}  
  \Phi^{(n)}_{k\ell}
  &=
  \sum_{i,j}
  (-1)^{\ell}
  \binom{j+\ell}{j}
  S^{(n)}_{ij}\gnbarijk{n}{i}{k}{j+\ell}.
\end{align*}

\begin{figure}[htbp]
  \centering
  \includegraphics{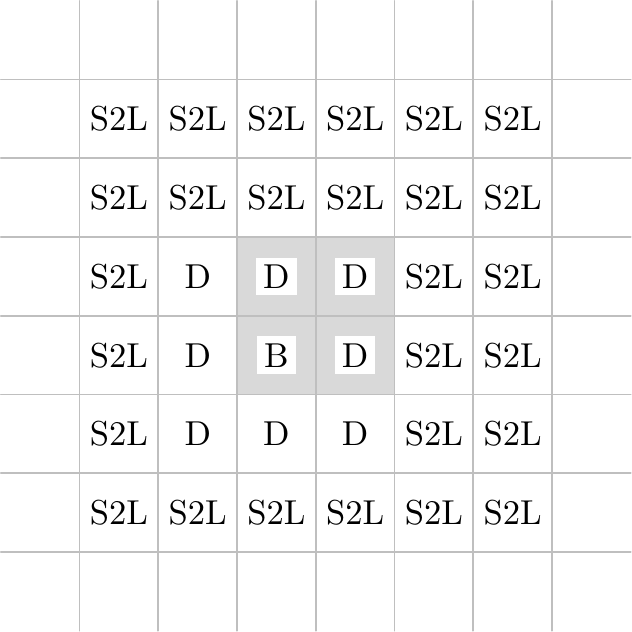}
  \caption{Interaction lists for box B. Parent box of B is shown
    shaded; boxes marked D have interactions evaluated directly from
    \cref{equ:analysis:sum}; boxes marked S2L have interactions
    evaluated via source-to-local operation. Unmarked boxes interact
    with B through its parent box so that the interactions need not be
    explicitly evaluated.}
  \label{fig:fmm:lists}
\end{figure}

In order to apply the shift operations, we enumerate candidate source
boxes which may contribute to the field in a box,
\cref{fig:fmm:lists}. This gives rise to two interaction lists, the D
list containing boxes which contribute via direct evaluation of the
field for each source and field point, and the S2L list, whose
contributions are evaluated using the S2L operation acting on source
and local expansion coefficients. \cref{fig:fmm:lists} indicates that
the D list is made up of neighbors of box B, leaving~27 boxes which
may contribute to the field in B via S2L operations. Contributions
from all other boxes are transferred into B from its parent box during
the downward pass. Using the axial translation invariance and the
symmetry in $r$ and $r_{1}$ reduces the number of Green's function
expansions to evaluated to twelve, those for boxes at $z\geq z_{1}$
and $r\geq r_{1}$. The expansions are identical for any value of
$z-z_{1}$ but must be updated for each $r_{1}$ during the downward
pass. Once generated, the expansion is used to update the local
expansion on the outer boxes, and to include the contribution of those
boxes' source terms to the local expansion on B.

Finally, we note that the modal Green's function can be written using
scaled coordinates, so that
\begin{align}
  \label{equ:gfunc:scaled}
  G^{(n)}(\sigma r,\sigma r_{1}, \sigma x)
  &=
  \frac{1}{\sigma}G^{(n)}(r, r_{1}, x),
\end{align}
which would allow for the shift operators to be precomputed and
generated at each level as required. This has been implemented but
found not to give a time saving, since each level of the tree requires
twice as many shift operators as its parent level, half of which are
new. In practice, we find that the bottleneck in the code is the
evaluation of direct interactions rather than the computation of the
S2L operators.

\section{Algorithm}
\label{sec:algorithm}

Combining the elements of the previous sections, we present an
algorithm for a uniform Fast Multipole Method in an axisymmetric
domain. Input is a list of $N_{s}$ source points $(r_{i},z_{i})$,
$i=1,\ldots,N_{s}$ and modal amplitudes $S_{i}^{(n)}$, $n=0,\ldots,N$,
and a list of $N_{f}$ field points $(r_{j},z_{j})$,
$j=1,\ldots,N_{f}$.

\begin{algorithm}
  \caption{Fast Multipole Method for cylindrical coordinate systems}
  \label{alg:fmm}
  \begin{algorithmic}
    \STATE{set tree depth $d$}

    \STATE{sort source points $(r_{i},z_{i})$ and field points
      $(r_{j},z_{j})$ by Morton index and assign to leaf nodes at
      level $d$}

    \STATE{calculate leaf box moments from source
      amplitudes~\cref{equ:analysis:moments}}\\

    \COMMENT{upward pass}

    \FOR{$\ell=d-1,\ldots,2$}     

    \STATE{evaluate box moments at level $\ell$ from child box moments
      at level $\ell+1$, \cref{equ:analysis:shift}}

    \ENDFOR\\

    \COMMENT{downward pass}

    \FOR{$\ell=2,\ldots,d$}

    \FOR{$i=0,\ldots,4^{\ell-1}-1$}

    \STATE{shift parent box local expansion from level $\ell-1$ to
      level $\ell$ boxes}

    \ENDFOR

    \FOR{$i=0,\ldots,2^{\ell}-1$}

    \STATE{evaluate coefficients of Green's function expansions for
      radial station 
      $i$}
    
    \FOR{$j=0,\ldots,2^{\ell}-1$}

    \STATE{apply S2L operators for FO and BO translations of source in
      box $(i,j)$}

    \STATE{apply S2L operators for FI and BI translation to update
      local expansion in box $(i,j)$}

    \ENDFOR
    \ENDFOR
    \ENDFOR

    \FOR{$i=0,\ldots,4^{d}-1$}

    \STATE{for field points in box $i$ evaluate sum of local expansion
      and direct contribution from neighbor boxes}
    
    \ENDFOR
  \end{algorithmic}
\end{algorithm}

The implementation in \cref{alg:fmm} is for a uniform FMM which does
not generate an adaptive decomposition when assigning points to boxes.
This was decided upon to reduce the number of shift operators
required in evaluating box interactions. 

\section{Results}
\label{sec:results}

\begin{figure}[htbp]
  \centering
  \includegraphics{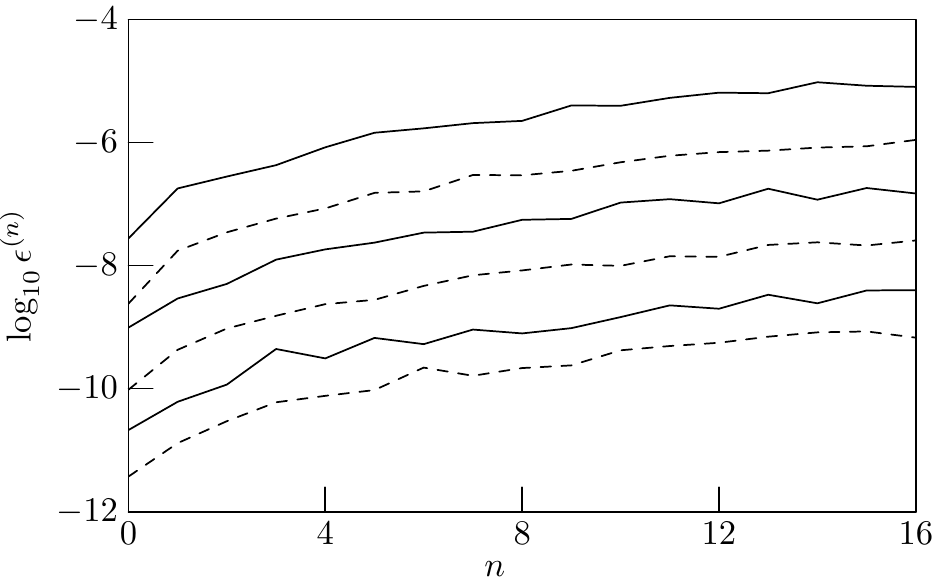}
  \caption{Error $\epsilon$ against mode order for real part of modal
    amplitude $N_{s}=2^{16}$, tree depth $d=6$, expansion order
    $M=6,8,10,12,14,16$ from top to bottom of plot.}
  \label{fig:results:error}
\end{figure}

\begin{figure}[htbp]
  \centering
  \includegraphics{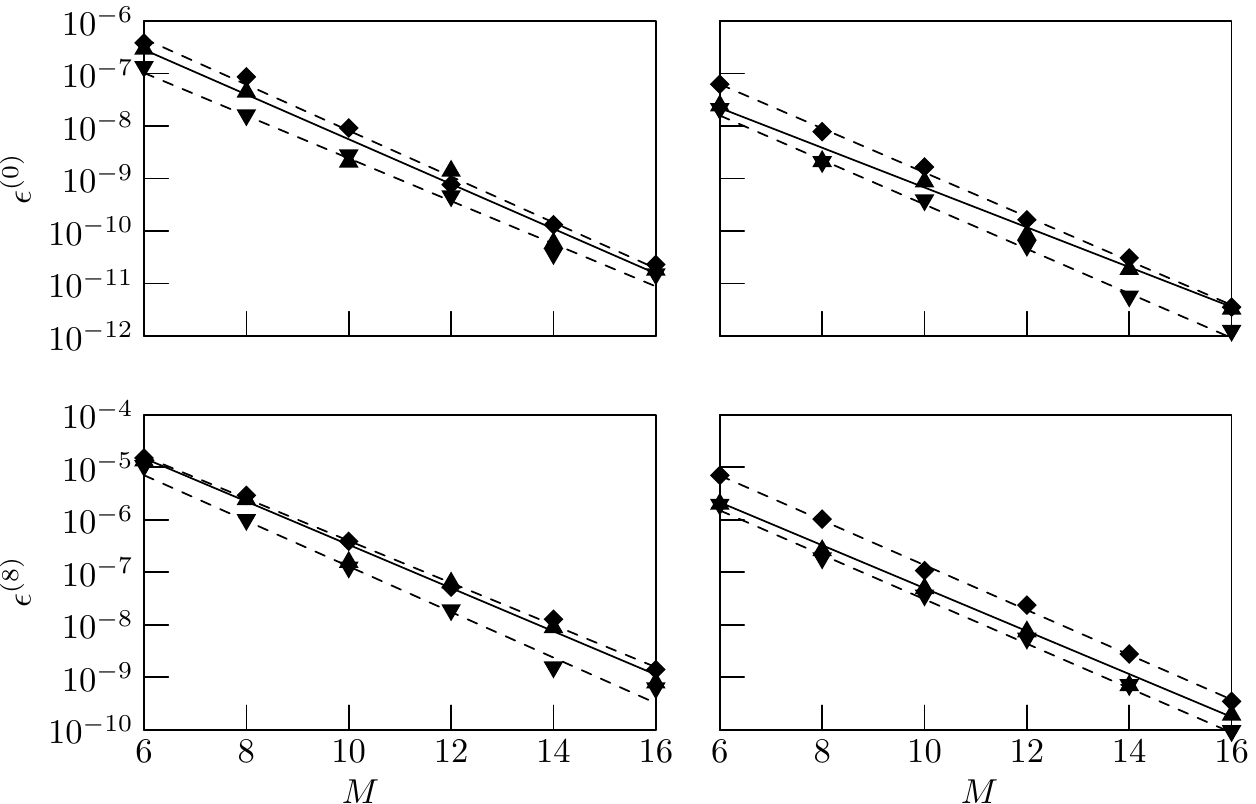}
  \caption{Error $\epsilon$ for $n=0$ (top row) and $n=8$ (bottom row)
    against expansion order $M$ for $N_{s}=2^{14}$ (left column) and
    $N_{s}=2^{16}$ (right column): diamonds, tree depth $d=5$; upward
    triangles, $d=6$; downward triangles, $d=7$; fitted lines
    $\epsilon\approx C^{-0.4M}$.}
  \label{fig:results:error:scaling}
\end{figure}

The algorithm has been tested for accuracy and computation time using
source and field points randomly distributed over $0\leq r,z\leq 1$,
with random modal amplitudes $0<S^{(n)}<1$, $n=0,\ldots,17$. In each
case, source number $N_{s}$ is set equal to number of field points
$N_{f}$, with $N_{s}=2^{q}$, $q=10,\ldots,16$, and the same order of
expansion is used for source and field terms. Results are presented
for varying $N_{s}$, maximum expansion order $M$, and tree depth
$d$. Code is written in GNU C, with \texttt{gcc} optimization
\texttt{-O3} and Goto BLAS matrix-vector operations. Calculations were
performed on one core of an Intel i5-6200U laptop running
at~2.3GHz. Similar code and optimizations were used for the direct
evaluations used as an error reference.

Error is evaluated for each modal amplitude of the field,
\begin{align}
  \label{equ:results:error}
  \epsilon^{(n)}
  &=
  \frac{\max|\Phi^{(n)}_{FMM}-\Phi^{(n)}_{D}|}{\max|\Phi_{D}^{(n)}|},
\end{align}
where $\Phi^{(n)}_{FMM}$ is modal amplitude evaluated using the new
algorithm, and $\Phi^{(n)}_{D}$ is that found by direct
evaluation. Sample results for error as a function of mode number and
expansion order are shown for $N_{s}=2^{16}$ in
\cref{fig:results:error}. The method is clearly accurate, especially
for higher order expansions, where eleven digit accuracy is achieved
for the axisymmetric mode. The error increases at larger $n$, where
the absolute value of the modal amplitudes is smaller, making the
relative error measure larger.

\cref{fig:results:error:scaling} shows the variation in error with
expansion order for the axisymmetric $n=0$ mode and for $n=8$. The
error scales approximately as $C^{-0.4M}$ with weak dependence on
tree depth. The algorithm performs well with respect to convergence
over the range of problem sizes tested here.

\begin{figure}[htbp]
  \centering
  \includegraphics{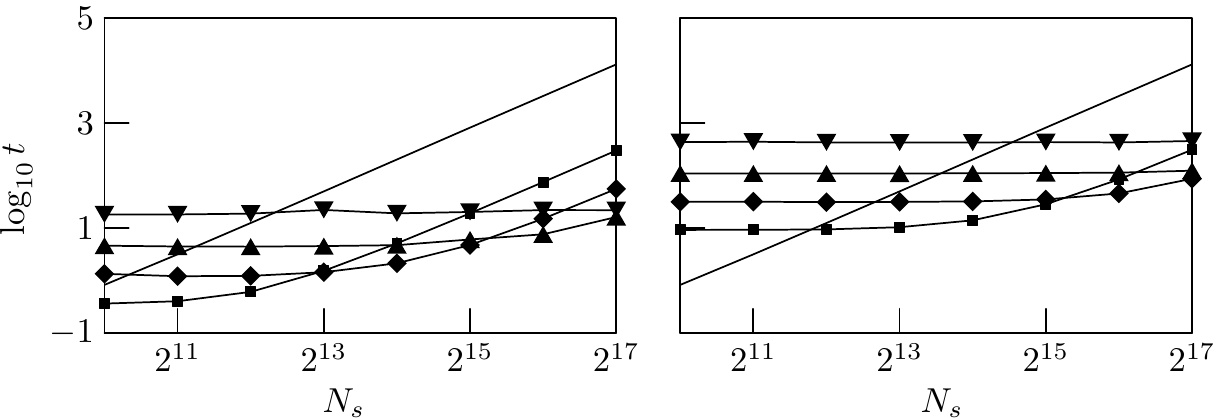}
  \caption{Computation time against source and field point number, for
    direct evaluation (solid line), and trees of depth~4 (boxes)
    depth~5 (diamonds),~6 (upward triangles), and~7 (downward
    triangles); left hand plot: expansion order $M=6$; right hand plot
    order expansion order $M=16$.}
  \label{fig:results:time}
\end{figure}

\cref{fig:results:time} shows basic data for computation time as a
function of problem size. The time for direct evaluation is shown and
scales at the expected $N_{f}N_{s}=N_{s}^{2}$ rate. The computation
time for the FMM algorithm behaves similarly for the low $M=6$ order
and high $M=16$ order cases, with times being shifted up by the change
in expansion order. The computation time in each case is roughly
constant for small $N_{s}$, where the evaluation time is dominated by
the set up cost, which depends on the tree depth. As the problem size
increases, the evaluation time for the downward pass begins to
dominate the calculation time which increases proportional to
$N_{s}^{2}$, but with a much smaller leading constant than for direct
evaluation. Again, this is the expected behavior as the direct
evaluation of near-field interactions becomes the largest part of the
calculation. With increasing tree depth, the box to box evaluations
become correspondingly faster as the number of sources per box becomes
smaller. 

\begin{figure}[htbp]
  \centering
  \includegraphics{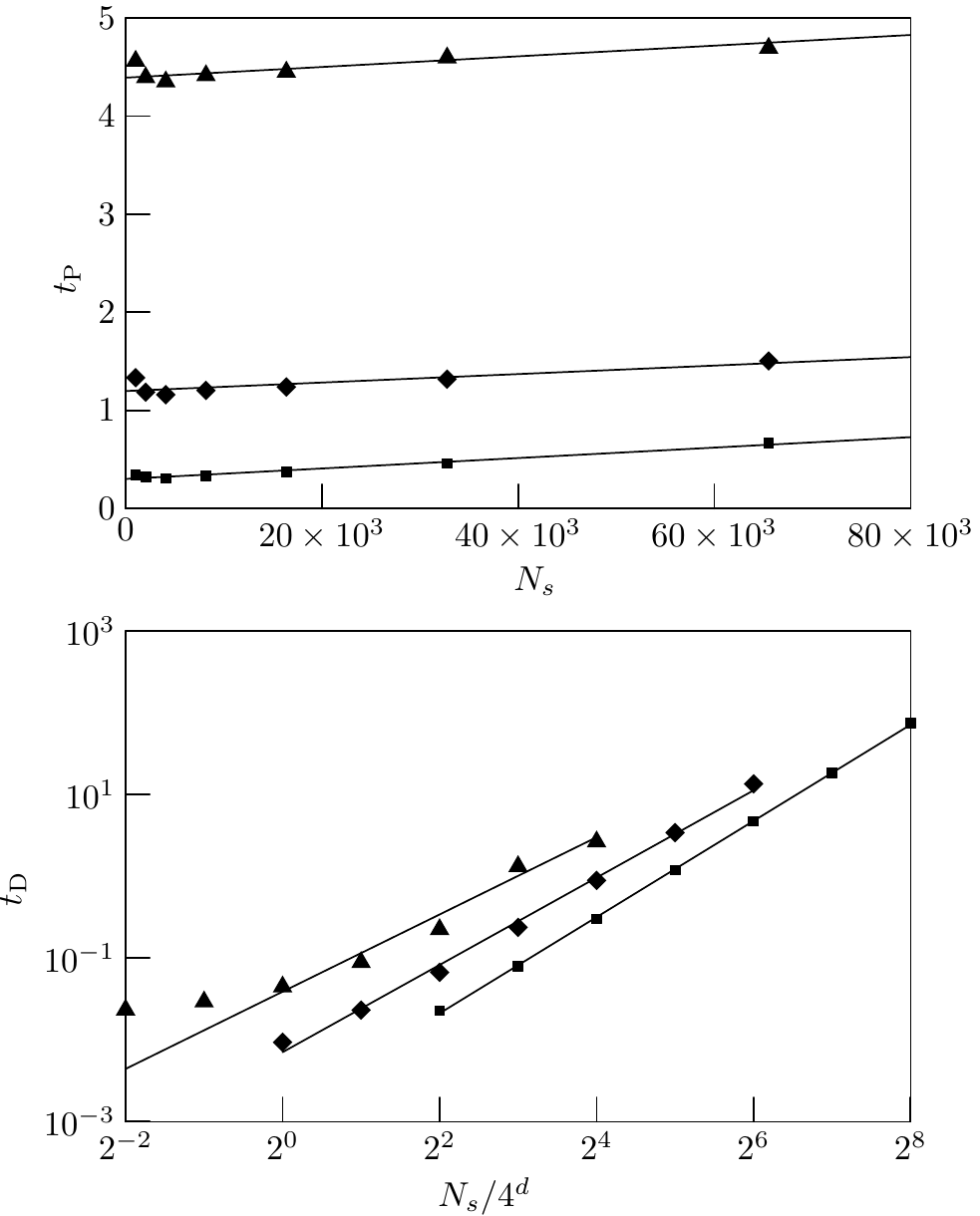}
  \caption{Execution time for phases of calculation, expansion order
    $M=6$; boxes: tree depth $d=4$; diamonds, $d=5$, triangles,
    $d=6$. Upper plot, time for upward and downward passes, with
    linear fit; lower plot, time for local field evaluation, power law
    fit to points with $N_{s}/4^{d}\geq1$.}
  \label{fig:results:time:breakdown}
\end{figure} 

\cref{fig:results:time:breakdown} shows the breakdown of computation
time between the two parts of the calculation, as a function of tree
depth. The initialization phase, made up of the upward and downward
passes, scales approximately linearly with problem size, with a
leading constant determined by the tree depth. Initialization time
increases with tree depth, as the number of boxes increases. The time
for local field evaluation, in the lower plot, scales well on
$N_{s}/4^{d}$, the average number of sources per box, with the time
reducing with tree depth. The implication is that computation time and
accuracy are determined by the balance between initialization and
local field evaluation, which depends on source number, expansion
order, and tree depth.

\begin{figure}[htbp]
  \centering
  \includegraphics{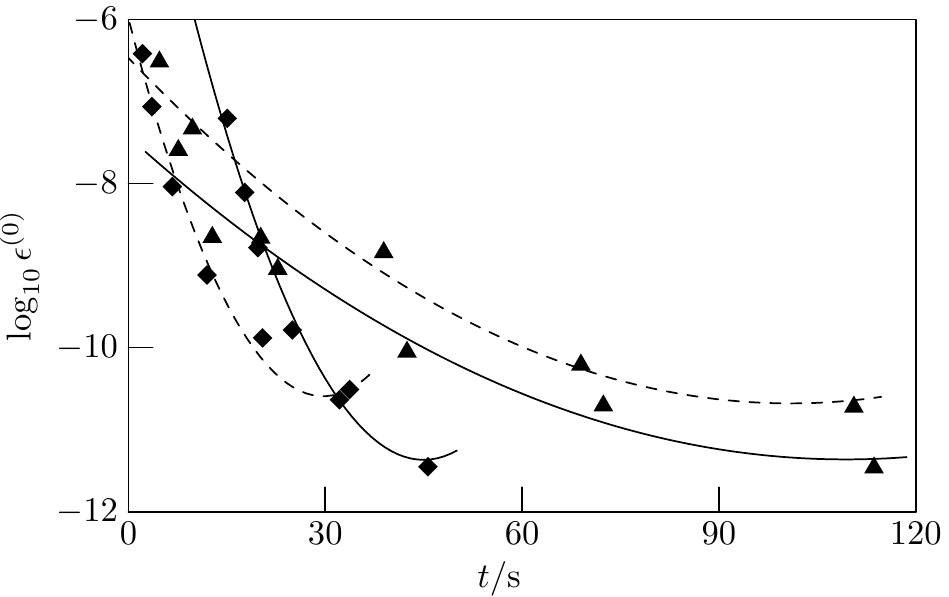}
  \caption{Error $\epsilon^{(0)}$ against computation time: diamonds
    depth $d=5$, triangles $d=6$; solid lines, $N_{s}=2^{16}$; dashed
    lines, $N_{s}=2^{14}$; curves are second order fits to
    $\log\epsilon$.}
  \label{fig:results:time:error}
\end{figure} 

\cref{fig:results:time:error} gives results illustrating this
balance, plotting error against computation time, found by varying
problem size and expansion order at two tree depths. For each tree
depth and source number, there is a trend towards a minimum error,
with the greater tree depth requiring a greater total computation time
for higher order accuracy. 


\section{Conclusions}
\label{sec:conclusions}

The Fast Multipole Method has been extended to non-axisymmetric
problems in cylindrical domains, by evaluating the amplitudes of the
modes in a Fourier expansion of the source and field in a Poisson
problem. Testing by comparison with direct evaluation shows
convergence to up to ten digit accuracy, and orders of magnitude
speed-up, depending on expansion order. Open questions remain. The
first is the efficient evaluation of the Legendre functions used to
find the modal Green's function, which is common to the direct and
fast methods, and constitutes the largest computational demand in the
method. A second is the formulation of the method in a form which
allows the use of BLAS level~3 operations, which should allow the code
to be optimized further. Finally, we note that the approach taken here
should be applicable to the Helmholtz problem, though with some
greater difficulty in evaluating the modal Green's functions. 




\appendix

\section{Evaluation of Green's functions and derivatives}
\label{sec:gfunc}

The modal Green's function $G^{(n)}(r,r_{1},x)$ is defined:
\begin{align}
  \label{equ:gfunc:def}
  G^{(n)}(r,r_{1},x)
  &=
  \int_{0}^{2\pi}
  \frac{\E^{\I n\theta_{1}}}{4\pi R}\,\D\theta_{1},\\
  R^{2} &= r^{2} + r_{1}^{2} - 2rr_{1}\cos\theta_{1} + x^{2}.\nonumber
\end{align}
Cohl and Tohline~\cite{cohl-tohline99} give an expansion for $1/R$,
\begin{align}
  \label{equ:gfunc:expansion}
  \frac{1}{R}
  &=
  \frac{1}{\pi\sqrt{rr_{1}}}
  \sum_{m=-\infty}^{\infty}
  \E^{\I m(\theta-\theta_{1})}
  Q_{m-1/2}(\chi),\\
  \chi &= \frac{r^{2}+r_{1}^{2}+x^{2}}{2rr_{1}},\nonumber
\end{align}
where $Q_{\nu}(\chi)$ is the Legendre function of the second
kind. Integration over $\theta_{1}$ yields
\begin{align}
  \label{equ:gfunc:basic}
  G^{(n)}(r,r_{1},x)
  &=
  \frac{Q_{n-1/2}(\chi)}{2\pi\sqrt{rr_{1}}}.   
\end{align}

Using the recursion for the Legendre
function~\cite[8.732.2]{gradshteyn-ryzhik80}, with the functional
dependence on coordinates suppressed for clarity,
\begin{align}
  \label{equ:gfunc:recursion}
  (2n-3)G^{(n-2)} &= 4(n-1)\chi G^{(n-1)} - (2n-1)G^{(n)}.
\end{align}
For $\chi>1$, the forward recursion is unstable and the backward
recursion is stable, but computationally
expensive~\cite{helsing-karlsson14}. To generate the sequence of modal
Green's functions, we apply the approach of Helsing and
Karlsson~\cite{helsing-karlsson14} and use the forward recursion for
$\chi<1.008$, beginning with the initial values~\cite{cohl-tohline99},
\begin{align}
  \label{equ:gfunc:Q0}
  Q_{-1/2}(\chi) &= \mu K(\mu),\\
  Q_{1/2}(\chi) &= \chi\mu K(\mu) - (1+\chi)\mu E(\mu),\\
  \mu &= \sqrt{\frac{2}{1+\chi}}.\nonumber
\end{align}
Here $K(\cdot)$ and $E(\cdot)$ are the complete elliptic integrals of
the first and second kind respectively. These are computed using the
method of Carlson~\cite{carlson95}.

For $\chi\geq1.008$, the backward recursion is used starting with
arbitrary values of $Q_{n-1/2}(\chi)$ and $Q_{n-3/2}(\chi)$ for
$n=N+80$, performing the downward recursion to $n=0$ and scaling the
sequence using the known value of $Q_{-1/2}(\chi)$,
\cref{equ:gfunc:Q0}. Values of $G^{(n)}$ evaluated using this
procedure have been checked against numerical integration and have
been found to be correct to machine precision.

Given values of $G^{(n)}$, $n=0,\ldots,N$, the derivatives of $G^{(n)}$
can be found using a combination of the recursion relations for the
Legendre function and the Laplace equation. For concision, we
introduce the notation
\begin{align}
  \label{equ:gfunc:gnijk}
  \gnijk{n}{i}{j}{k}
  &=
  \frac{\partial^{i+j+k}}{\partial r^{i}\partial r_{1}^{j}\partial
    x^{k}} G^{(n)}(r,r_{1},x),\\
  \label{equ:gfunc:gnbar}
  \gnbarijk{n}{i}{j}{k}
  &=
  \frac{1}{i!}\frac{1}{j!}\frac{1}{k!}
  \frac{\partial^{i+j+k}}{\partial r^{i}\partial r_{1}^{j}\partial
    x^{k}} G^{(n)}(r,r_{1},x),
\end{align}
so that the Taylor series for $G^{(n)}$ is given by
\begin{align}
  \label{equ:gfunc:taylor}
  G^{(n)}(r+\Delta r,r_{1}+\Delta r_{1},x+\Delta x)
  &=
  \sum_{m=0}^{\infty}
  \sum_{i+j+k=m}
  \gnbarijk{n}{i}{j}{k}(\Delta r)^{i}(\Delta r_{1})^{j}(\Delta
  x)^{k}.
\end{align}

The derivatives are evaluated using a recursion based on the Laplace
equation, similar to the approach of Strickland and Amos who used the
axisymmetric stream function
equation~\cite{strickland-amos90,strickland-amos92}. Here we use the
Laplace equation for a field with azimuthal dependence
$\exp \I n\theta$. This recursion requires starting values which can
be found using the properties of the Legendre
functions~\cite[8.732]{gradshteyn-ryzhik80}. For derivatives with
respect to $x$,
\begin{align}
  \label{equ:gfunc:derivatives:x}
  \gnbarijk{n}{0}{0}{1}
  &=
  (n-1/2)
  \left[
    \left(\gnbarijk{n}{0}{0}{0} + \gnbarijk{n-1}{0}{0}{0}\right)
    \frac{x}{\rho_{+}^{2}}
    +
    \left(\gnbarijk{n}{0}{0}{0} - \gnbarijk{n-1}{0}{0}{0}\right)
    \frac{x}{\rho_{-}^{2}}
  \right],\\
  \gnbarijk{n}{0}{0}{k+1}
  &=
  (n-1/2)
  \sum_{q=0}^{k}
  \frac{1}{q!(k+1)}
  \begin{aligned}[t]
    \biggl[
    \left(\gnbarijk{n}{0}{0}{k-q} + \gnbarijk{n-1}{0}{0}{k-q}\right)
      \frac{\partial^{q}}{\partial x^{q}}
      \left(
        \frac{x}{\rho_{+}^{2}}
      \right)
      +\\
      \left(\gnbarijk{n}{0}{0}{k-q} - \gnbarijk{n-1}{0}{0}{k-q}\right)
      \frac{\partial^{q}}{\partial x^{q}}
      \left(
        \frac{x}{\rho_{-}^{2}}
      \right)
    \biggr],
    \end{aligned}\\
  \rho_{\pm}^{2} &= (r\pm r_{1})^{2} + x^{2}.\nonumber
\end{align}
To evaluate derivatives for $n=0$, the relation
$\gnijk{-1}{i}{j}{k}\equiv\gnijk{1}{i}{j}{k}$ can be used.

For the derivatives with respect to $r$,
\begin{align}
  \gnbarijk{n}{1}{0}{0}
  &=
  -\frac{\gnbarijk{n}{0}{0}{0}}{2r}
  + (n-1/2)\frac{r^{2} - r_{1}^{2} - x^{2}}{2r}
  \begin{aligned}[t]
    \biggl[
    \left(\gnbarijk{n}{0}{0}{0} + \gnbarijk{n-1}{0}{0}{0}\right)
    \frac{1}{\rho_{+}^{2}}
    +\\
    \left(\gnbarijk{n}{0}{0}{0} - \gnbarijk{n-1}{0}{0}{0}\right)
    \frac{1}{\rho_{-}^{2}}
    \biggr],    
  \end{aligned}\\
  \gnbarijk{n}{1}{0}{k}
  &=
  -\frac{\gnbarijk{n}{0}{0}{k}}{2r}
  + \frac{n-1/2}{2r}
  \sum_{q=0}^{k}
  \frac{1}{q!}
  \begin{aligned}[t]
     \bigg[
    &\left(\gnbarijk{n}{0}{0}{k-q} + \gnbarijk{n-1}{0}{0}{k-q}\right)
    \frac{\partial^{q}}{\partial x^{q}}
    \frac{r^{2} - r_{1}^{2} - x^{2}}{\rho_{+}^{2}} + \\
    &\left(\gnbarijk{n}{0}{0}{k-q} - \gnbarijk{n-1}{0}{0}{k-q}\right)
    \frac{\partial^{q}}{\partial x^{q}}
    \frac{r^{2} - r_{1}^{2} - x^{2}}{\rho_{-}^{2}}    
    \bigg].
  \end{aligned}
\end{align}
Derivatives $\gnijk{n}{0}{1}{k}$ are found by exchanging $r$ and
$r_{1}$, with a corresponding swap of indices in the derivatives.

Finally,
\begin{align}
  \gnbarijk{n}{1}{1}{k}
  &=
  -\frac{\gnbarijk{n}{0}{1}{k}}{2r}\\
  &+\frac{n-1/2}{2r}
  \sum_{q=0}^{k}
  \frac{1}{q!}
  \begin{aligned}[t]
    \biggl[
    &
    \left(\gnbarijk{n}{0}{0}{k-q} + \gnbarijk{n-1}{0}{0}{k-q}\right)
    \frac{\partial^{q+1}}{\partial r_{1}\partial x^{q}}
    \frac{r^{2} - r_{1}^{2} - x^{2}}{\rho_{+}^{2}}
    + \\
    &
    \left(\gnbarijk{n}{0}{1}{k-q} + \gnbarijk{n-1}{0}{1}{k-q}\right)
    \frac{\partial^{q}}{\partial x^{q}}
    \frac{r^{2} - r_{1}^{2} - x^{2}}{\rho_{+}^{2}}
    + \\
    &
    \left(\gnbarijk{n}{0}{0}{k-q} - \gnbarijk{n-1}{0}{0}{k-q}\right)
    \frac{\partial^{q+1}}{\partial r_{1}\partial x^{q}}
    \frac{r^{2} - r_{1}^{2} - x^{2}}{\rho_{-}^{2}}
    + \\
    &
    \left(\gnbarijk{n}{0}{1}{k-q} - \gnbarijk{n-1}{0}{1}{k-q}\right)
    \frac{\partial^{q}}{\partial x^{q}}
    \frac{r^{2} - r_{1}^{2} - x^{2}}{\rho_{-}^{2}}
    \biggr].\nonumber
  \end{aligned}  
\end{align}

To compute the remaining derivatives, we make use of a recursion based
on the Laplace equation in cylindrical coordinates. Noting that
$\partial/\partial\theta\to\I n$,
\begin{align}
  \label{equ:gfunc:laplace}
  \frac{1}{r}
  \frac{\partial}{\partial r}
  \left(
    r
    \frac{\partial G^{(n)}}{\partial r}
  \right)
  -
  \frac{n^{2}}{r^{2}}G^{(n)}
  +
  \frac{\partial^{2}G^{(n)}}{\partial x^{2}} = 0.
\end{align}

This yields the relation,
\begin{align}
  \label{equ:gfunc:laplace:recursion}
  \gnbarijk{n}{i+2}{j}{k}
  &=
  \sum_{u=0}^{i}
  \frac{(-1)^{u}}{r^{u+1}}
  \frac{1}{(i+1)(i+2)}
  \left[
    \frac{(u+1)n^{2}}{r}\gnbarijk{n}{i-u}{j}{k}
    -
    (i-u+1)\gnbarijk{n}{i-u+1}{j}{k}    
  \right],
\end{align}
which can be used to recursively generate higher derivatives with
respect to $r$. Switching $r$ and $r_{1}$ gives a corresponding
relation for the higher derivatives with respect to $r_{1}$ and allows
a complete set of derivatives to be evaluated to any required order.



\end{document}